\documentclass[12pt]{article}
      \usepackage{latexsym,amsfonts,amssymb}
      \usepackage{tikz-cd}
      \usepackage{amsmath, esint}
      \newcommand{\overbar}[1]{\mkern 2.5mu\overline{\mkern-2.5mu#1\mkern-2.5mu}\mkern 2.5mu}      
      \newtheorem{theorem}{Theorem}[section]

      \newtheorem{lemma}[theorem]{Lemma}

      \newtheorem{Ex}{Example}[section]
      \newtheorem{remark}{Remark}[section]
      
      \def\nn{\nonumber}
      \def\rf#1{\mbox{$(\ref{#1})$}}

      \def\be{\begin{equation}} 
      \def\ee{\end{equation}} 
      \def\beqn{\begin{eqnarray}} 
      \def\eeqn{\end{eqnarray}} 
      \def\beq{\begin{eqnarray*}} 
      \def\eeq{\end{eqnarray*}}
      \def\proof{{\bf Proof:}\ }
      \def\ep{\varepsilon} 
      \def\la{\lambda} 
      \def\ra{\rightarrow} 
      \def\sta{\stackrel} 

\begin{document}
      \title{\bf Large Deviations For Randomly Weighted Sums of Random Measures\thanks{Research supported by
      the Natural Science and Engineering Research Council of Canada}}
      \author{Shui Feng\\Department of Mathematics and Statistics\\McMaster
      University\\Hamilton, Ontario
      Canada L8S 4K1\\
      shuifeng@mcmaster.ca}
      \date{\empty}
      \maketitle
      \begin{abstract}
     
      Let $\{{\bf \mathcal{Z}}_n:n\geq 1\}$ be a sequence of i.i.d. random probability measures. Independently, for each $n\geq 1$, let  $(X_{n1},\ldots, X_{nn})$ be a random vector of positive random variables that add up to one. This paper studies the large deviation principles for the randomly weighted sum $\sum_{i=1}^{n} X_{ni} \mathcal{Z}_i$.  In the case of finite Dirichlet weighted sum of Dirac measures, we obtain an explicit form for the rate function.
     It  provides a new measurement of divergence between probabilities.   As applications, we obtain the large deviation principles for a class of randomly weighted means including the Dirichlet mean and the corresponding posterior mean.  We also identify the minima of relative entropy with mean constraint in both forward and reverse directions.

      \vspace*{.125in} \noindent {\bf Key words:} Bayesian non-parametrics, Dirichlet process, Dirichlet means, gamma process, large deviation, relative entropy.
      \vspace*{.125in}

      \noindent {\bf MSC2020-Mathematics Subject Classifications:}
      Primary: 60G57; Secondary: 62F15.

      \end{abstract}
      \section{Introduction}
      \setcounter{equation}{0}

      Let $(\Omega,{\cal F}, P)$ be a probability space, $S$ be a compact subset of real line $\mathbb{R}$, and $\mathcal{S}$ be  the Borel $\sigma$-algebra. Let $C(S)$ and $B(S)$ denote the spaces of  continuous functions and bounded measurable functions on $S$, respectively. Let $M_1(S)$ denote the  space of probability measures on $(S, \mathcal{S})$ equipped with the weak topology.

      For any $n \geq 1$, set
      \[
      \triangle_n=\{(x_1, \ldots,x_n): x_i\geq 0, \sum_{k=1}^n x_k=1\}.
      \]
         
      Let $\{\mathcal{Z}_i:i\geq 1\}$ be i.i.d. $M_1(S)$-valued  random probability measures with common distribution $\Pi$, and, independently, ${\bf W}:=\{(W_{n1},W_{n2}, \ldots, W_{nn}): n \geq 1\}$ be a triangular array of positive random variables. Set
      
\[
X_{ni} =\frac{W_{ni}}{\sum_{k=1}^n W_{nk}},\  i=1,\ldots,n.
\]
and

   \be\label{rm-2}
      \mathcal{W}_n({\Pi};{\bf W})=\sum_{i=1}^n X_{ni}\mathcal{Z}_i.    \ee 
      
\vspace{2mm}

      The random distributions $\mathcal{Z}_i$ usually carry the information of location distributions and the coefficients $X_{ni}$ are the corresponding weights or proportions.  The collection of random measures $\mathcal{W}_n({\Pi},{\bf W})$ is a subclass in the family of randomly weighted averages. The latter has drawn a great deal of interests from researchers  due to a wide range of applications including Bayesian nonparametrics, genetics, and  the path behaviour of diffusion processes.   See the survey papers \cite{Lijoi-Prunster09},\cite{Pitman2018}, and the references therein for a comprehensive coverage on the historical and recent development of the subjects.

      The main concern of this paper is the large deviations of $ \mathcal{W}_n({\Pi};{\bf W})$ as $n$ tends to infinity.    To motivate the study, consider the following subclass of models.  
       
       For each $\theta >0$, let ${\bf W}$ be the triangular array where for each $n\geq 1$ $W_1^n, \ldots, W_n^n$ are i.i.d. exponential with parameter $\theta/n$.  Independently,  let $\xi_1,\xi_2,\ldots$ be i.i.d. with common distribution $\nu_0$ in $M_1(S)$, $\delta_x$ be the Dirac measure concentrated at $x$, and $\mathcal{Z}_i =\delta_{\xi_i}$.  

Then, as $\theta$ tends to infinity, $\mathcal{W}_n(\Pi, {\bf W})$ converges in distribution to the empirical distribution of $\xi_1, \ldots, \xi_n$

 \[
    \mathcal{L}_{n,\nu_0}= \frac{1}{n}\sum_{i=1}^n \delta_{\xi_i}.
      \] 
      
\vspace{2mm}

Now letting $n$ go to infinity, we obtain the limit $\nu_0$.       
       
 On the other hand,  let $U_1,U_2,\ldots$ be i.i.d. with Beta$(1,\theta)$ distribution.   The Dirichlet process (\cite{Fer73}) with base distribution $\nu_0$ and concentration parameter $\theta$ is the random measure
 
      \be \label{DIR}
     \mathcal{V}_{\theta,\nu_0}= \sum_{i=1}^{\infty}V_i\delta_{\xi_i},
      \ee
      where $V_1=U_1, V_i =(1-U_1)\cdots (1-U_{i-1})U_i,  i\geq 2.$
  
It is known (\cite{ish-zare02}) that $\mathcal{W}_n(\Pi, {\bf W})$ converges in distribution to the Dirichlet process $\mathcal{V}_{\theta,\nu_0}$ as $n$ tends to infinity.  Taking the limit of $\theta$ going to infinity,  we obtain the same limit $\nu_0$.

Instead of taking limits of $\theta$ and $n$ separately, we can consider the special case $\theta=n$.  Denote $\mathcal{W}_n(\Pi,{\bf W})$ in this case by

\be\label{expo-weight}
\mathcal{W}_{n,\nu_0}=\sum_{i=1}^n X_{ni}\delta_{\xi_i}.
\ee
 
 \vspace{2mm}
 
 One can show that $\mathcal{W}_{n,\nu_0}$ also converges in distribution to $\nu_0$ as $n$ tends to infinity.  The following diagram is a summery of these relations.

\vspace{0.4cm}
 \be\label{dia1}
 \begin{tikzcd}[row sep=14pt,column sep={18mm,between origins}]
&\mathcal{W}_n(\Pi,{\bf W})  \arrow[dddd,"\begin{matrix}\theta\vspace{-1mm}\\ \vspace{-1mm}\big\downarrow \\\vspace{-2mm} \infty\end{matrix}"] \arrow[ddddrrr, "\theta=n\ra \infty"]  \arrow[rrr,"n\ra\infty"] & &  & \mathcal{V}_{\theta,\nu_0}\arrow[dddd,"\begin{matrix}\theta\vspace{-1mm}\\ \vspace{-1mm}\big\downarrow \\\vspace{-2mm} \infty\end{matrix}"] & \\
 &  & &  &  & \\
    &  & &  &  & \\[4mm]
 &  & &  &  & \\
    & \mathcal{L}_{n,\nu_0} \arrow[rrr,"n\ra \infty"]& &  & \nu_0& \\
\end{tikzcd}
\ee 
  \vspace{0.4cm}    
     
    It is natural to investigate and to compare the asymptotic behaviours of the three different estimation of  $\nu_0$.  The large deviations of $\mathcal{L}_{n,\nu_0}$ from $\nu_0$ when $n$ tends to infinity is given by Sanov's theorem (cf. \cite{Sanov57}) with the following good rate function (\cite{DV75I})
   
    \be\label{re-1}
    H(\mu|\nu_0)= \sup\{\langle\mu,f\rangle-\log \langle \nu_0,e^f\rangle: f\in C(S)\},   
    \ee
        where $\langle\mu,f\rangle= \int f(x)\mu(d\,x)$.   The function $H(\mu|\nu_0)$ is the relative entropy of $\mu$ with respect to $\nu_0$.

         When $\theta$ converges to infinity, the large deviations of $\mathcal{V}_{\theta,\nu_0}$ from $\nu_0$ (\cite{DF01},\cite{LS87})  has a good rate function given by the reversed relative entropy $H(\nu_0|\mu)$ when the support of $\mu$ is contained in the support of $\nu_0$.
     
    Our new result of large deviations for $\mathcal{W}_{n,\nu_0}$ corresponds to the diagonal limit in \rf{dia1}. The good rate function will be identified as the minimum of the sum of relative entropies with respect to $\mu$ and $\nu_0$. Like the relative entropy, it provides an alternate measurement of divergence between probability measures. Comparisons among the three good rate functions reveal  refined differences of the three sequences of random measures.   
      
   The random probability measure $\mathcal{W}_{n,\nu_0}$ also appears in the Dirichlet posterior. More specifically, let $\eta_1, \eta_2 ,\ldots, $ be i.i.d. given the Dirichlet process $\mathcal{V}_{\theta,\nu_0}$. Then the posterior distribution of $\mathcal{V}_{\theta,\nu_0}$ given $\eta_1, \ldots,\eta_n$ is
   
\[
U \mathcal{V}_{\theta,\nu_0}+(1-U)\sum_{i=1}^n X_{i}\delta_{\eta_i} 
\]
where $U$ is Beta$(\theta,n)$ distributed, $(X_{n1},\ldots,X_{nn})$ has Dirichlet$(1,\ldots,1)$ distribution, and all random variables appearing on the right hand side are independent given $\eta_1, \ldots,\eta_n$.  Clearly the term $\sum_{i=1}^n X_{ni}\delta_{\eta_i} $ is the same as  $\mathcal{W}_{n,\nu_0}$ if $\eta_1, \ldots,\eta_n$ are i.i.d. with common distribution $\nu_0$.  If $\theta$ is fixed, then the large deviation principle for $\mathcal{W}_{n,\nu_0}$ is associated with the annealed large deviation principle for the Dirichlet posterior.

 The development of the paper is as follows. We starts with the large deviation principles for $ \mathcal{W}_n({\Pi};{\bf W})$ in Section 2. Our focus will then shift in Section 3 to the study of the large deviations of $\mathcal{W}_{n,\nu_0}$, and its connections to Sanov's theorem and the large deviations of  $\mathcal{V}_{\theta,\nu_0}$.  Several applications will be discussed in Section 4 including large deviations for Dirichlet and Dirichlet posterior means.  We also obtain some minimization properties of the relative entropy.   The logarithm $\log$ in this article is the natural logarithm, with base $e$.    All terms and definitions regarding large deviations are found in the reference \cite{DeZe98}.

\section{Large Deviations for  $ \mathcal{W}_n({\Pi};{\bf W})$}

   \setcounter{equation}{0}  
   
 In this section, we discuss the large deviations for 
     $ \mathcal{W}_n({\Pi};{\bf W})$ defined in \rf{rm-2} with $S=[0,1]$.   We also assume in the sequel that
     
     \be\label{assum-1}
     \lim_{n\ra \infty}\mathbb{E}\bigg[\sum_{i=1}^n X^2_{ni}\bigg]=0.
     \ee     

\vspace{2mm}

Let $\{f_i \in C(S):  i=1, 2,\ldots\}$ be a countable dense subset of $C(S)$. The metric on $M_1(S)$ is  given by

\[
\rho(\upsilon,\mu)=\sum_{i=1}^{\infty}\frac{|\langle\upsilon-\mu, f_i\rangle|\wedge 1}{2^i}. 
\] 

\vspace{2mm}

Let $\nu_0$ denote the mean measure of $\mathcal{Z}_1$. Then for any $f$ in $C(S)$

\[
\mathbb{E}\bigg[\langle \mathcal{W}_n({\Pi};{\bf W}), f\rangle\bigg] =\langle\nu_0,f\rangle,\ \ 
\mbox{Var}\bigg[\langle \mathcal{W}_n({\Pi};{\bf W}), f\rangle\bigg] = \mathbb{E}[\sum_{i=1}^n X^2_{ni}] \mbox{Var}\bigg(\langle \mathcal{Z}_1,f\rangle\bigg).
\]

\vspace{2mm}

It follows from \rf{assum-1} that for any $\ep>0$

\[
\lim_{n\ra \infty}P\{ \rho(\mathcal{W}_n({\Pi};{\bf W}), \nu_0)\geq \ep\}=0.
\]

\vspace{2mm}

Thus $ \mathcal{W}_n({\Pi};{\bf W})$ converges in probability to $\nu_0$ in $M_1(S)$. For any $\delta >0$ and any $ \mu$ in $M_1(S)$, let

\[
B(\mu,\delta)=\{\upsilon \in M_1(S): \rho(\upsilon,\mu)< \delta\}, \  \ \overbar{B}(\mu,\delta)=\{\upsilon \in M_1(S): \rho(\upsilon,\mu)\leq \delta\}, 
\]
and 
\[
S_{\mu}=\{t\in S: t \ \mbox{is a continuity point of}\ \mu\}.
\]
\vspace{2mm}

For any $m\geq 1$, and any $0=t_0<t_1<\ldots<t_m<1$  in  $S_{\mu}$, set

\beq
&&B_{t_1,\ldots,t_m}(\mu,\delta)=\{\upsilon \in M_1(S): \rho_{t_1,\ldots, t_m}(\upsilon,\mu)< \delta\}\\
&& \overbar{B}_{t_1,\ldots, t_m}(\mu,\delta)=\{\upsilon \in M_1(S): \rho_{t_1,\ldots, t_m}(\upsilon,\mu)\leq \delta\}, 
\eeq
where 

 \[
 \rho_{t_1,\ldots, t_m}(\upsilon,\mu)=\upsilon([t_{m},1])-\mu([t_{m},1])|+\sum_{i=1}^{m} |\upsilon([t_{i-1},t_i))-\mu([t_{i-1},t_i))|.
 \]
 \vspace{2mm}

 Then the following holds.

\begin{theorem}\label{grm}
For each $n \geq 1$, let $\Pi_n$ denote the law of $\mathcal{W}_n(\Pi, {\bf W})$. Assume that for any $\mu $ in $M_1(S)$,  $m\geq 1$, and any  $0<t_1<\ldots<t_m<1$  in $S_{\mu}$, there exists $I_{t_1,\ldots,t_m}(\mu)$ such that   
\beqn\label{grm-eq1}
 &&\lim_{\delta \ra 0}\liminf_{n\ra \infty}\frac{1}{n}\log \Pi_n\{  B_{t_1,\ldots,t_m}(\mu,\delta)  \}\nn\\
 &&= \lim_{\delta \ra 0}\limsup_{n\ra \infty}\frac{1}{n}\log \Pi_n\{   \overbar{B}_{t_1,\ldots,t_m}(\mu,\delta) \}\\
 &&=-I_{t_1,\ldots, t_m}(\mu). \nn
 \eeqn
 Then the sequence $\{\Pi_n:n\geq 1\}$ satisfies a large deviation principle with good rate function 
 \[
 I(\mu)=\sup_{t_1,\ldots, t_m \in S_{\mu}}\{I_{t_1,\ldots, t_m}(\mu)\}.
 \]
 \end{theorem}
\proof Since the space $M_1(S)$ is compact, the sequence $\{\Pi_n: n \geq 1\}$ is exponentially tight.   It follows from Theorem (P) in \cite{Puk91} that  every subsequence of $\{\Pi_n:n \geq 1\}$ has a further subsequence that satisfies a large deviation principle with a good rate function.  To prove the theorem it suffices to prove that all these good rate functions are the same as $I(\cdot)$.   This is clearly true if we can show that  

\beq
 &&\lim_{\delta \ra 0}\liminf_{n\ra \infty}\frac{1}{n}\log \Pi_n\{ B(\mu,\delta)\}\nn\\
 &&= \lim_{\delta \ra 0}\limsup_{n\ra \infty}\frac{1}{n}\log \Pi_n\{\overbar{B}(\mu,\delta) \}\\
 &&=-I(\mu). \nn
 \eeq

 For any $\mu$ in $M_1(S)$ there exist $f_1, \cdots, f_r $ in $C(S)$ and $\ep >0$ such that
\[
\{\upsilon \in M_1(S):|\langle \upsilon, f_j \rangle -\langle \mu, f_j \rangle|< \ep:j=1,\cdots, r \} \subset B(\mu,\delta).
\]

Let 
\[
c = \sup\{|f_j(x)|: x \in S, j =1,\cdots,r\},
\]
and choose the continuity points $0<t_1<\cdots< t_m<1$ in  $S_{\mu}$ such that
\[
 \sup\{|f_j(x)-f_j(y)|: x,y \in [t_{i-1},t_{i}],i=1,\ldots, m+1; j=1,\cdots,r \}< \ep/4
\]
where $t_0=0, t_{m+1}=1$.
Let $0<\delta_1 <\frac{\ep}{4c}$. Then 
for any $\upsilon$ in $B_{t_1, \cdots, t_m}(\mu, \delta_1)$  and any $f_j$, we have
\beq
|\langle \upsilon, f_j \rangle -\langle \mu, f_j \rangle|&=&|\int_{[t_{m},1]}[f_j(x)-f_j(t_{m})+ f_j(t_{m})](\upsilon(dx)-\mu(dx))\\
&  +& \hspace{0.3cm}\sum_{i=1}^{m}\int_{[t_{i-1},t_{i})}[f_j(x)-f_j(t_{i-1})+f_j(t_{i-1})](\upsilon(dx)-\mu(dx))|\\
&<& \frac{\ep}{2} +  c\delta_1  <  \ep,
\eeq
which implies that 
\[
 B_{t_1, \cdots, t_m}(\mu, \delta_1)\subset
\{\mu \in M_1(S):|\langle \upsilon, f_j \rangle -\langle \mu, f_j \rangle|< \ep:j=1,\cdots, r \} \subset B(\mu,\delta).
\]

Noting that the continuity points $t_1, \ldots, t_m$ does not depend on $\delta_1$, it follows from \rf{grm-eq1}  that
\beq
-I(\mu)&\leq& -I_{t_1,\ldots,t_m}(\mu)\\
&=&\lim_{\delta_1 \ra 0}\liminf_{n\ra \infty}\frac{1}{n}\log \Pi_n\{  B_{t_1,\ldots,t_m}(\mu,\delta_1)\} \\
 &\leq&  \liminf_{n\ra \infty}\frac{1}{n}\log 
\Pi_n\{ B(\mu,\delta)\}.\eeq

It now follows by letting $\delta$ go to zero that
\be \label{local-1}
\liminf_{\delta\ra 0}\liminf_{n\ra \infty}\frac{1}{n}\log 
\Pi_n\{ B(\mu,\delta)\}\geq-I(\mu).\ee

On the other hand, for any $\mu$ in $M_1(S)$, $\delta_2>0$, $m\geq 1$, and any continuity points $0<t_1<\cdots<t_m<1$ of $\mu$, the set $B_{t_1,\ldots,t_m}(\mu,\delta_2)$ is open in $M_1(S)$. One can choose $\delta$ small enough so that
\[
B(\mu, \delta)\subset \overbar{B}(\mu,\delta)\subset B_{t_1,\ldots,t_m}(\mu,\delta_2)\subset \overbar{B}_{t_1,\ldots,t_m}(\mu,\delta_2).\]

Letting $\delta$ go to zero we obtain that
\[
\limsup_{\delta\ra 0}\limsup_{n\ra \infty}\frac{1}{n}\log 
\Pi_n\{ \overbar{B}(\mu,\delta)\}\leq \limsup_{n\ra \infty}\frac{1}{n}\log 
\Pi_n\{ \overbar{B}_{t_1,\ldots,t_m}(\mu,\delta_2)\}
\]

Letting $\delta_2$ go to zero followed by taking the infimum of $-I_{t_1,\ldots,t_m}(\mu)$ over  the continuity points of $\mu$ we obtain 
\be\label{local-2}
\limsup_{\delta\ra 0}\limsup_{n\ra \infty}\frac{1}{n}\log 
\Pi_n\{ \overbar{B}(\mu,\delta)\}\leq -I(\mu).
\ee

Putting together \rf{local-1} and \rf{local-2} we obtain the result.

\hfill $\Box$

\begin{remark}\label{re-2.1}
The space $S$ may be generalized to a compact set in a Polish space.  But the compactness is crucial in general.   
\end{remark}

\begin{remark}\label{re-2.2}
The limit $I_{t_1,\ldots,t_m}(\mu)$ appearing in \rf{grm-eq1} is just a non-negative number. There is no other conditions imposed.  This makes the verification a bit easier. 
\end{remark}


Theorem~\ref{grm} works well for random probabilities that have explicit forms when projected to finite partitions of $S$.  We demonstrate this through two examples. 

\begin{Ex}\label{ex-empirical}
Let  $\xi_1,\ldots, \xi_n,\ldots$ be i.i.d. with common law $\nu_0$, $W_{n1}=\cdots=W_{nn}=1$, and $\mathcal{Z}_i=\delta_{\xi_i}$. Then $X_{n1}=\cdots=X_{nn}=\frac{1}{n}$, $\Pi$ is the law of $\delta_{\xi_1}$, and  $\mathcal{W}_n(\Pi, {\bf W})$ is the usual empirical distribution
\[
\mathcal{L}_{n,\nu_0}=\frac{1}{n}\sum_{i=1}^n\delta_{\xi_i}.
\]
Assume that $\nu_0$ is positive on all open intervals. For any $m\geq 2$ and any $$0=t_0<t_1<\cdots<t_{m-1}<1,$$ let $$p_m=\nu_0([t_{m-1},1]),\  p_i=\nu_0([t_{i-1}, t_i)), i=1, \ldots, m-1.$$  We have $$\bigg(\mathcal{W}_n(\Pi, {\bf W})([0,t_1)),\ldots,\mathcal{W}_n(\Pi, {\bf W})([t_{m-1},1])\bigg)\sta{d}{=}\frac{1}{n}(Y_1^n, \ldots, Y_m^n)$$ where  $\sta{d}{=}$ denotes equality in distribution, and $(Y_1^n, \ldots, Y_m^n)$ is multinomial with parameters $n, p_1, \ldots, p_m$. Then the condition \rf{grm-eq1} holds for $$I_{t_1, \ldots,t_m} (\mu) =\sum_{i=1}^m q_i\ln \frac{q_i}{p_i}$$ where 
$$q_m =\mu([t_{m-1}, 1]),\ q_i=\mu([t_{i-1},t_i), i=1, \dots, m-1.$$

\end{Ex}

\begin{Ex}\label{ex-Dirichlet}
Let $\gamma(t)$ denote the gamma subordinator on $[0,+\infty)$ and set
\[
W_i = \gamma(i)-\gamma(i-1), i=1, 2, \ldots.
\]
 Let $\{J_{ik}: k\geq 1\}$ denote all the jump sizes of $\gamma(t)$ over the interval $[i-1,i]$. Then it follows from the Gamma-Dirichlet algebra \emph{(Theorem~1.1 in \cite{Feng10})} that the Dirichlet process, defined in \rf{DIR},  with base distribution $\nu_0$ on $[0,1]$ and concentration parameter $n$ has the form
 \[
\mathcal{V}_{n,\nu_0}= \frac{1}{\gamma(n)}\sum_{i=1}^n\sum_{k=1}^{\infty}J_{ik}\delta_{\xi_{ik}}
 \]
 where $\{\xi_{ik}:i=1, \ldots,n; k=1,\ldots\}$ are i.i.d. with common distribution $\nu_0$. Reorganizing the terms and applying the Gamma-Dirichlet algebra again we obtain
\beqn\label{decom}
\mathcal{V}_{n,\nu_0}&=&\frac{1}{\sum_{l=1}^n W_l}\sum_{i=1}^n W_i\sum_{k=1}^\infty\frac{J_{ik}}{W_i}\delta_{\xi_{ik}}\\
&=&\frac{1}{\sum_{l=1}^n W_l}\sum_{i=1}^n W_i\mathcal{Z}_i\nn\eeqn 
 where $\mathcal{Z}_1, \ldots, \mathcal{Z}_n$ are i.i.d. Dirchlet processes with base distribution $\nu_0$ and concentration parameter one. Thus the Dirichlet process $\mathcal{V}_{n,\nu_0}$ corresponds to $\mathcal{W}_n(\Pi, {\bf W})$ where $W_1, W_2, \ldots$ are i.i.d. exponential with parameter one and each $\mathcal{Z}_i$ is a copy of $\mathcal{V}_{1,\nu_0}$. Assume that the support of $\nu_0$ is $[0,1]$. Then the condition \rf{grm-eq1} holds \emph{(\cite{DF01})} for
 
 \[
I_{t_1, \ldots,t_m} (\mu)=\sum_{i=1}^m p_i\ln\frac{p_i}{q_i}
 \] 
 where $p_i$'s and $q_i$'s are defined as in Example~\ref{ex-empirical}. Note that the order between $p_i$'s and $q_i$'s is switched here.  
\end{Ex}

\section{Large Deviations for $\mathcal{W}_{n,\nu_0}$ }

   \setcounter{equation}{0}  
   
In this section we establish the large deviations for $\mathcal{W}_{n,\nu_0}$ defined in \rf{expo-weight}. The key in the proof  is to verify the condition \rf{grm-eq1} in Theorem~\ref{grm}.  Noting that the weights of $\mathcal{W}_{n,\nu_0}$ are the same as the weights of the Dirichlet process $\mathcal{V}_{n,\nu_0}$ in \rf{decom}. But the location distributions are different. This innocent-looking difference leads to very different large deviation behaviours as shown below.

For any $m\geq 1$ and any partition $A_1, \ldots, A_m$ of $S$, define the map
\[
\pi: M_1(S)\ra \triangle_m, \upsilon \ra (\upsilon(A_1), \ldots, \upsilon(A_m)).
\]

Let $p_i =\nu_0(A_i), i=1,\ldots, m$, and 
\[
n_i =\#\{1\leq k \leq n: \xi_k \in A_i\}, \  i=1, \ldots, m.
\]

  Then $\pi(\mathcal{W}_{n,\nu_0})$ has Dirichlet$(n_1,\ldots,n_m)$ distribution, and  $\pi(\mathcal{L}_{n,\nu_0})$ has multinomial distribution with parameters $n$, $p_1, \ldots,p_m$. Without loss of generality, we assume  $n_i\geq 1, p_i>0$ for $i=1,\ldots,m$. For  $\mu, \upsilon$ in $M_1(S)$, set $$\pi(\mu)=(q_1,\ldots,q_n),\pi(\upsilon)=(o_1,\ldots,o_n).$$  Define
  \[
  H(\pi(\upsilon)|\pi(\mu))=\sum_{i=1}^m o_i\log\frac{o_i}{q_i},\  H(\pi(\upsilon)|\pi(\nu_0))=\sum_{i=1}^m o_i\log\frac{o_i}{p_i},\]
   where $0\log0 =0\log\frac{0}{0}=0$.   

\begin{lemma}
For any  $\mu, \upsilon$ in $M_1(S)$, we have 
\beqn\ \label{gen-local}
&&\lim_{\delta \ra 0}\liminf_{n\ra \infty}\frac{1}{n}\log P\bigg\{\bigg|\bigg(\pi(\mathcal{W}_{n,\nu_0}), \pi(\mathcal{L}_{n,\nu_0})\bigg)- \bigg(\pi(\mu),\pi(\upsilon)\bigg)\bigg|<\delta\bigg\}\nn\\
&&=\lim_{\delta \ra 0}\limsup_{n\ra \infty}\frac{1}{n}\log P\bigg\{\bigg|\bigg(\pi(\mathcal{W}_{n,\nu_0}), \pi(\mathcal{L}_{n,\nu_0})\bigg)- \bigg(\pi(\mu),\pi(\upsilon)\bigg)\bigg|\leq \delta\}
\\
&&=-\bigg[H(\pi(\upsilon)|\pi(\mu))+H(\pi(\upsilon)|\pi(\nu_0))\bigg]\nn
\eeqn
where
\[
\bigg|\bigg(\pi(\mathcal{W}_{n,\nu_0}), \pi(\mathcal{L}_{n,\nu_0})\bigg)- \bigg(\pi(\mu),\pi(\upsilon)\bigg)\bigg|=\sum_{i=1}^m\bigg[\bigg|\mathcal{W}_{n,\nu_0}(A_i)-q_i\bigg|+\bigg|\mathcal{L}_{n,\nu_0}(A_i)-o_i\bigg|\bigg].
\]

\end{lemma}
\proof We only give the proof of equality

\beqn
&&\lim_{\delta \ra 0}\liminf_{n\ra \infty}\frac{1}{n}\log P\bigg\{\bigg|\bigg(\pi(\mathcal{W}_{n,\nu_0}), \pi(\mathcal{L}_{n,\nu_0})\bigg)- \bigg(\pi(\mu),\pi(\upsilon)\bigg)\bigg|<\delta\bigg\}\label{open-ball}\\
&&\hspace{2cm}=-\bigg[H(\pi(\upsilon)|\pi(\mu))+H(\pi(\upsilon)|\pi(\nu_0))\bigg].\nn
\eeqn

 The proof for $P\bigg\{\bigg|\bigg(\pi(\mathcal{W}_{n,\nu_0}), \pi(\mathcal{L}_{n,\nu_0})\bigg)- \bigg(\pi(\mu),\pi(\upsilon)\bigg)\bigg|\leq\delta\bigg\}$ is similar. 

For  any $\delta>0,$ $1\leq l\leq m$, and  $ r_i\geq 1, \sum_{i=1}^l r_i=n$, let
\[
F(r_1,\ldots, r_l;\delta)=  \frac{\Gamma(n)}{\Gamma(r_1)\cdots\Gamma(r_l)}\idotsint\limits_{{\bf D}_{\delta, l}}\prod_{i=1}^l x_i^{r_i-1}d\,x_1\ldots d\,x_{l-1}
\]
where
\[
{\bf D}_{\delta, l}=\{(x_1,\ldots, x_l)\in \triangle_l:\sum_{i=1}^l|x_i-q_i|<\delta\}.
\]

We can extend the domain of $F$ to zero integers by defining
\[
F(n_1, \ldots,n_m;\delta)=F(n_{k_1}, \ldots, n_{k_r};\delta)
\] 
where $\{n_{k_j}:j=1, \ldots,r\}$ consists of all non-zero elements of $\{n_i:i=1,\ldots,m\}$.

For any $\mu$ in $M_1(S)$, we have  
\beqn
&&\hspace{-2.5cm}P\bigg\{\bigg|\bigg(\pi(\mathcal{W}_{n,\nu_0}), \pi(\mathcal{L}_{n,\nu_0})\bigg)- \bigg(\pi(\mu),\pi(\upsilon)\bigg)\bigg|<\delta\bigg\}\label{expan}\\
\nn\\
&&\hspace{3.5cm} =\sum_{\sum_{k=1}^m|\frac{n_i}{n}-o_i|<\delta}A(n_1,\ldots,n_m;\delta)\nn
\eeqn
 where
 \[
A(n_1,\ldots,n_m;\delta)={n \choose n_1 \cdots n_m} F(n_1,\ldots, n_m;\delta)\prod_{i=1}^m p_i^{n_i}. \]

 Noting that $A(n_1,\ldots,n_m;\delta)$ depends on positive $n_i$'s only. We thus assume, without loss of generality, that  $n_i\geq 1$ for $i=1,\ldots, m$.  The arguments for general cases are similar except  the replacement of  $(n_1, \ldots, n_m)$ with its subset consisting of non-zero elements. 

By Stirling's formula for the gamma function, there exist positive constants $c_1<c_2$ such that for all $n\geq 1$ and $ n_i\geq 1, \sum_{i=1}^m n_i=n$
  \[
  c_1 \prod_{i=1}^m\bigg(\frac{n_i}{n}\bigg)^{-n_i} \sqrt{\frac{n_1\ldots n_m}{n}}\leq \frac{\Gamma(n)}{\Gamma(n_1)\cdots\Gamma(n_m)} \leq c_2  \prod_{i=1}^m\bigg(\frac{n_i}{n}\bigg)^{-n_i} \sqrt{\frac{n_1\ldots n_m}{n}}.\]
  
 Noting that for any integer $n\geq 1$
 \[
 \sqrt{2\pi}n^{n +\frac{1}{2}}e^{-n} \leq n! \leq  e n^{n +\frac{1}{2}}e^{-n}, \]
 it follows that there exist constants $c_3<c_4$ such that
 \[
 c_3 \prod_{i=1}^m\bigg(\frac{n_i}{n}\bigg)^{-n_i} \sqrt{\frac{n}{n_1\ldots n_m}}\leq {n \choose n_1 \ldots n_m}\leq c_4 \prod_{i=1}^m\bigg(\frac{n_i}{n}\bigg)^{-n_i} \sqrt{\frac{n}{n_1\ldots n_m}}.\]

 Putting these together it follows that
 \[
 c_5 \prod_{i=1}^m\bigg(\frac{n_i}{n}\bigg)^{-2n_i}\leq {n \choose n_1 \ldots n_m}\frac{\Gamma(n)}{\Gamma(n_1)\cdots\Gamma(n_m)}\leq c_6 \prod_{i=1}^m\bigg(\frac{n_i}{n}\bigg)^{-2n_i} \]
for some positive constants $c_5<c_6$.

 For any ${\bf x}=(x_1,\ldots, x_m), {\bf y}=(y_1,\ldots, y_m)$ in $\triangle_m$, define the function
 \[
 \Psi({\bf x},{\bf y})= H({\bf y}|{\bf x})+H({\bf y}|\pi(\nu_0)).
 \]

 Then we have  for  ${\bf y}=\left(\frac{n_1}{n}, \ldots, \frac{n_m}{n}\right)$ 
 \beqn
&& c_5\idotsint\limits_{{\bf D}_{\delta, m}}\exp\{-n(\Psi({\bf x}, {\bf y})+n^{-1}\log \prod_{i=1}^m x_i)\}d\,x_1,\cdots d\,x_{m-1} \nn\\
&&\nn\\
&& \ \ \ \ \hspace{1.8cm} \leq  A(n_1,\ldots,n_m;\delta)\label{upper-lower}\\
&&\nn\\
&&\ \ \ \ \leq  c_6\idotsint\limits_{{\bf D}_{\delta, m}}\exp\{-n(\Psi({\bf x}, {\bf y})+n^{-1}\log \prod_{i=1}^m x_i\}d\,x_1,\cdots d\,x_{m-1}\nn\eeqn
where
\[
{\bf D}_{\delta, m}=\{(x_1,\ldots, x_m)\in \triangle_m:\sum_{i=1}^m |x_i-q_i|<\delta\}.\]
  
  We now divide the upper estimations into three cases.
  
  {\bf Case 1}. $q_i>0$ for all $i$.  For any $\ep>0$ one can choose $n$ large and $\delta$ small such that
  \[
 \bigg| \Psi({\bf x}, {\bf y})+n^{-1}\log \prod_{i=1}^m x_i-\Psi(\pi(\mu), \pi(\upsilon))-n^{-1}\log \prod_{i=1}^m q_i\bigg|\leq \ep.  \]
  
 It follows that 
  \beqn\label{gen-local1}
  &&A(n_1,\ldots,n_m;\delta)\nn\\
  && \hspace{2cm}\leq   c_6 [2\delta]^m \exp\bigg\{-n\bigg(\Psi(\pi(\mu), \pi(\upsilon))+n^{-1}\log \prod_{i=1}^m q_i- \ep\bigg)\bigg\}. \eeqn
  
Since the total number of terms in \rf{expan} is at most $(n+1)^m$ and $\ep$ is arbitrary, we obtain that 
\beqn\label{gen-local2}
&&\lim_{\delta \ra 0}\limsup_{n\ra \infty}\frac{1}{n}\log P\bigg\{\bigg|\bigg(\pi(\mathcal{W}_{n,\nu_0}), \pi(\mathcal{L}_{n,\nu_0})\bigg)- \bigg(\pi(\mu),\pi(\upsilon)\bigg)\bigg|<\delta\bigg\}\nn\\
&& \hspace{2cm}\leq -\sum_{i=1}^lo_i\bigg(\log\frac{o_i}{q_i}+\log\frac{o_i}{p_i}\bigg) \\
&& \hspace{2cm}=- \bigg[H(\pi(\upsilon)|\pi(\mu))+H(\pi(\upsilon)|\pi(\nu_0))\bigg].\nn
\eeqn

 Next we turn to the situation where $q_i =0$ for some $i$. Without loss of generality, we assume that there exists an $1<l<m$ such that $q_i> 0$ for $1\leq i \leq l$  and $q_i=0$ for $i>l$. 
 
 {\bf Case 2}. $o_i >0$ for some $i>l$. In this case the term $\frac{n_i-1}{n}\log\frac{n_i/n}{x_i}$ in $$\Psi({\bf x}, {\bf y})+n^{-1}\log \prod_{i=1}^m x_i$$ converges to infinity  as $\delta$ tends to zero. Thus  
 \[
 \lim_{\delta \ra 0}\limsup_{n\ra \infty}\frac{1}{n}\log P\bigg\{\bigg|\bigg(\pi(\mathcal{W}_{n,\nu_0}), \pi(\mathcal{L}_{n,\nu_0})\bigg)- \bigg(\pi(\mu),\pi(\upsilon)\bigg)\bigg|\leq \delta\}\leq -\infty. \]

The result \rf{open-ball} follows from the fact that

\[
H(\pi(\upsilon)|\pi(\mu))+H(\pi(\upsilon)|\pi(\nu_0))=\infty,\]

{\bf Case 3}. $o_i=0$ for $i>l$.  It follows from direct calculation that
\beq
&&\exp\bigg\{-n\bigg(\Psi({\bf x}, {\bf y})+n^{-1}\log\prod_{i=1}^m x_i\bigg)\bigg\}\\
&&=\bigg(\prod_{i=l+1}^m x_i^{n_i-1}\bigg)\exp\bigg\{-n\bigg(\sum_{i=1}^ly_i\bigg(\log\frac{y_i}{x_i}+\log\frac{y_i}{p_i}\bigg)+n^{-1}\log \prod_{i=1}^l x_i\bigg)\bigg\}\\
&&\hspace{2cm} \times \exp\bigg\{-n\bigg(\sum_{i=l+1}^m y_i\bigg(\log y_i+\log\frac{y_i}{p_i}\bigg)\bigg)\bigg\}\\
&&\leq \exp\bigg\{-n\bigg(\sum_{i=1}^ly_i\bigg(\log\frac{y_i}{x_i}+\log\frac{y_i}{p_i}\bigg)+n^{-1}\log \prod_{i=1}^l x_i\bigg)\bigg\}\\
&& \hspace{2cm}\times \exp\bigg\{-n\bigg(\sum_{i=l+1}^m y_i\bigg(\log y_i+\log\frac{y_i}{p_i}\bigg)\bigg)\bigg\}
\eeq
on the domain ${\bf D}_{\delta,m}$. The exponential term
\[
\exp\bigg\{-n\bigg(\sum_{i=1}^ly_i\bigg(\log\frac{y_i}{x_i}+\log\frac{y_i}{p_i}\bigg)+n^{-1}\log \prod_{i=1}^l x_i\bigg)\bigg\}\]
can be estimated by an argument similar that used in deriving \rf{gen-local1}. 

For the second exponential term we have 

\beq
&&\lim_{\delta\ra 0}\limsup_{n\ra \infty}\frac{1}{n}\log  \exp\bigg\{-n\bigg(\sum_{i=l+1}^m y_i\bigg(\log y_i+\log\frac{y_i}{p_i}\bigg)\bigg)\bigg\}\\
&&\ \ \ \ \leq -\lim_{\delta \ra 0}\inf_{\sum_{i=l+1}^m y_i\leq\delta }\bigg\{\sum_{i=l+1}^m y_i\bigg(\log y_i+\log\frac{y_i}{p_i}\bigg)\bigg\}\\
&& \ \ \ \ =0.
\eeq

Thus \rf{gen-local2} also holds in this case.  It remains to check the lower bound in {\bf Cases } {\bf 1}  and  {\bf 3}. 
\vspace{2mm}

In {\bf Case 1}, the function $\Psi({\bf x}, {\bf y})+n^{-1}\log \prod_{i=1}^m x_i$ is continuous at $({\bf q}, {\bf o})$. Thus for any $\ep >0$ one can choose $\delta$ small so that
for 
\[
|({\bf x}, {\bf y})-({\bf q},{\bf o})|<\delta
\]
one  has 
\[
|\Psi({\bf x}, {\bf y})-\Psi({\bf q}, {\bf o}) + n^{-1}(\log\prod_{i=1}^m x_i-\log\prod_{i=1}^mq_i)|< \ep\]

By \rf{upper-lower} we have 
\beq
&&A(n_1,\ldots,n_m;\delta) \\
&& \geq c_5\idotsint\limits_{{\bf D}_{\delta,m}}exp\{-n[\Psi({\bf x}, {\bf y})+n^{-1}\log \prod_{i=1}^m x_i]\}d\,x_1,\cdots d\,x_m \\
&& \geq   c_5 e^{-n\ep}\exp\{ -n[\Psi({\bf q}, {\bf o})+n^{-1}\prod_{i=1}^m q_i] \} \idotsint\limits_{{\bf D}_{\delta,m}}d\,x_1\cdots d\,x_m
\eeq
which implies

\beqn\label{gen-lower-local2}
&&\lim_{\delta \ra 0}\liminf_{n\ra \infty}\frac{1}{n}\log P\bigg\{\bigg|\bigg(\pi(\mathcal{W}_{n,\nu_0}), \pi(\mathcal{L}_{n,\nu_0})\bigg)- \bigg(\pi(\mu),\pi(\upsilon)\bigg)\bigg|<\delta\bigg\}\\
&&\hspace{2cm} \geq- \bigg[H(\pi(\upsilon)|\pi(\mu))+H(\pi(\upsilon)|\pi(\nu_0))\bigg].\nn
\eeqn 

In {\bf Case 3},  let $T=\{1\leq i\leq m: q_i=0\} (o_i=0 \ \mbox{for}\ i \in T)$ and define
\[
\Psi_1({\bf x},{\bf y}) =\sum_{i \not\in T}y_i\log\frac{y_i}{x_i}, \ \Psi_2({\bf x},{\bf y}) =\sum_{i \in T}y_i\log\frac{y_i}{x_i}.\]

The function $\Psi_1({\bf x}, {\bf y})+n^{-1}\log\prod_{i\not\in T}x_i$ is clearly continuous at $({\bf q},{\bf o})$.  On the other hand, set  

\[
{\bf C}_{\delta,m}=\{{\bf x}=(x,\ldots,x_m)\in \triangle_m: \frac{1}{2}(q_i +\frac{\delta}{m})< x_i < q_i+\frac{\delta}{m}, \ \mbox{for all}\  i\}.
\]
Then the following holds on ${\bf C}_{\delta,m}$.
 
\[
\Psi_2({\bf x}, {\bf y})+n^{-1}\log\prod_{i\in T}x_i\leq  -m\delta\log\frac{\delta}{2m} +\frac{|T|}{n}\log\frac{\delta}{m} \]
where $|T|$ is the cardinality of $T$.

Since ${\bf C}_{\delta,m}$ is a subset of ${\bf D}_{\delta,m}$, it follows from \rf{upper-lower} and an argument similar to {\bf Case 1} that \rf{gen-lower-local2} holds in  {\bf Case 3}.

Putting together \rf{gen-local2} and \rf{gen-lower-local2} we obtain \rf{open-ball} and the lemma.

\hfill $\Box$

For any $\mu$ in $M_1(S)$,
and  any $0<t_1<\cdots<t_{m}<1$ in $S_{\mu}$, define
\[
\pi_{t_1,\ldots,t_m}(\mu)=(\mu([0,t_1)), \ldots, \mu([t_k,t_{k+1})), \ldots, \mu([t_{m},1])).
\]

\begin{lemma}\label{entropy-l1}
For any $\upsilon, \mu$ in $M_1(S)$,  we have
\beqn\label{entropy-e1}
&&H(\upsilon|\mu)+H(\upsilon|\nu_0)=\sup_{0<t_1<\cdots<t_m<1\  \in\ S_{\upsilon}\cap S_{\mu}}\{H(\pi_{t_1,\ldots,t_m}(\upsilon)|\pi_{t_1,\ldots,t_m}(\mu))\\
&&\hspace{5.3cm}+H(\pi_{t_1, \ldots,t_m}(\upsilon)|\pi_{t_1,\ldots,t_m}(\nu_0)\}.\nn
\eeqn
\end{lemma}

\proof It is known (\cite{Geo88}, \cite{DF01}) that for any $\upsilon, \mu$ in $M_1(S)$
\[
H(\upsilon|\mu) =\sup_{0<t_1<\cdots<t_m<1\  \in\  S_{\mu}}\{H(\pi_{t_1,\ldots,t_m}(\upsilon)|\pi_{t_1,\ldots,t_m}(\mu)\}.\]
Since  the supremum of sums is less than or equal to  the sum of supremums, it follows that
\beq
H(\upsilon|\mu)+H(\upsilon|\nu_0)
&\geq& \sup_{0<t_1<\cdots<t_m<1\  \mbox{in}\ S_{\upsilon} \cap S_{\mu}}\{H(\pi_{t_1,\ldots,t_m}(\upsilon)|\pi_{t_1,\ldots,t_m}(\mu))\\
&&\hspace{1.5cm}+H(\pi_{t_1, \ldots,t_m}(\upsilon)|\pi_{t_1,\ldots,t_m}(\nu_0)\}.\eeq

 To prove the other direction, we first recall the variational form \rf{re-1} of the relative entropy
\[
H(\upsilon|\mu)=\sup_{g \in C(S)}\{\langle\upsilon, g\rangle-\log \langle\mu, e^g\rangle\}.\]

For any $\ep >0$, there are $g,h$ in $C(S)$ such that
\[
H(\upsilon|\mu)\leq \langle\upsilon, g\rangle-\log \langle\mu, e^g\rangle+\ep
\]
and 
\[
H(\upsilon|\nu_0)\leq \langle\upsilon, h\rangle-\log \langle\mu, e^h\rangle+\ep.
\]

Since $S_{\upsilon}\cap S_{\mu}$ is dense in $S$, there exist $0<t_{n1}<\cdots<t_{nn}<1$ in $S_{\upsilon}\cap S_{\mu}$ such that
\beq
&&\lim_{n\ra \infty}\max_{1\leq i\leq n+1}\bigg\{\max\{|g(x)-g(y)|: x,y \in [t_{n(i-1)},t_{ni}]\bigg\}=0, \\
&&\lim_{n\ra \infty}\max_{1\leq i\leq n+1}\bigg\{\max\{|h(x)-h(y)|: x,y \in [t_{n(i-1)},t_{ni}]\}\bigg\}=0.
\eeq
where $t_{n0}=0, t_{n(n+1)}=1$. Set
\[
\alpha_{ni}=g(t_{ni}), \beta_{ni}=h(t_{ni}), i=1, \ldots, n+1
\]
and $$A_{n(n+1)}=[t_{nn},1], A_{ni}= [t_{n(i-1)},t_{ni}), i=1, \ldots, n.$$
Then there exists $c_n(g,h)$ such that
\beq
&&\lim_{n\ra\infty}c_n(g,h)=0\\
&&H(\upsilon|\mu)\leq \sum_{i=1}^{n+1}\alpha_{ni}\upsilon(A_{ni})-\log \sum_{i=1}^{n+1}e^{\alpha_{ni}}\mu(A_{ni})+\ep +c_n(g,h)\\
&&\hspace{1.4cm}\leq H(\pi_{t_{n1},\ldots,t_{nn}}(\upsilon)|\pi_{t_{n1},\ldots,t_{nn}}(\mu)) +\ep +c_n(g,h)
\eeq
and 
\beq
&&H(\upsilon|\nu_0)\leq \sum_{i=1}^{n+1}\beta_{ni}\upsilon(A_{ni})-\log \sum_{i=1}^{n+1}e^{\beta_{ni}}\nu_0(A_{ni})+\ep +c_n(g,h)\\
&&\hspace{1.6cm}\leq H(\pi_{t_{n1},\ldots,t_{nn}}(\upsilon)|\pi_{t_{n1},\ldots,t_{nn}}(\nu_0)) +\ep +c_n(g,h).
\eeq

Putting all these together we obtain \rf{entropy-e1}.

\hfill $\Box$

\begin{theorem}\label{DW}
Assume that the topological support of $\nu_0$ is $S$. Then the family of the laws of $\{\mathcal{W}_{n,\nu_0}:n\geq 1\}$ satisfies a large deviation principle on $M_1(S)$ with good rate function
\be\label{DW-rate}
J(\mu,\nu_0)=\inf_{\upsilon \in M_1(S)}\{H(\upsilon|\mu)+H(\upsilon|\nu_0)\}.
\ee
\end{theorem}

 \proof For any $\upsilon,\mu$ in $M_1(S)$, and any  $0<t_1<\cdots<t_{m}<1$ in $S_{\upsilon}\cap S_{\mu}$, let $$A_{m+1}=[t_{m}, 1], A_i=[t_{i-1}, t_{i}), i=1, \ldots, m.$$ Then it follows from \rf{gen-local} that
 
 \beq
 &&\lim_{\delta \ra 0}\liminf_{n\ra \infty}\frac{1}{n}\log P\bigg\{\bigg|\bigg(\pi_{t_1,\ldots,t_{m}}(\mathcal{W}_{n,\nu_0}), \pi_{t_1, \ldots,t_{m}}(\mathcal{L}_{n,\nu_0})\bigg)\\
 &&\hspace{5cm}- \bigg(\pi_{t_1,\ldots,t_{m}}(\mu),\pi_{t_1,\ldots,t_{m}}(\upsilon)\bigg)\bigg|<\delta\bigg\}\\
&&= \lim_{\delta \ra 0}\limsup_{n\ra \infty}\frac{1}{n}\log P\bigg\{\bigg|\bigg(\pi_{t_1,\ldots,t_{m}}(\mathcal{W}_{n,\nu_0}), \pi_{t_1, \ldots,t_{m}}(\mathcal{L}_{n,\nu_0})\bigg)\\
&&\hspace{5cm}- \bigg(\pi_{t_1,\ldots,t_{m}}(\mu),\pi_{t_1,\ldots,t_{m}}(\upsilon)\bigg)\bigg|\leq \delta\bigg\}
\\
&&=-\bigg[H\bigg(\pi_{t_1,\ldots, t_{m}}(\upsilon)|\pi_{t_1,\ldots, t_{m}}(\mu)\bigg)+H\bigg(\pi_{t_1,\ldots, t_{m}}(\upsilon)|\pi_{t_1,\ldots, t_{m}}(\nu_0)\bigg)\bigg].
 \eeq
 
 Since $M_1(S)\times M_1(S)$ is compact, it follows, by Lemma~\ref{entropy-l1} and an arguments similar to that used in the proof of Theorem~\ref{grm}, that 
 the family of the laws of  $(\mathcal{W}_{n,\nu_0},  \mathcal{L}_{n,\nu_0})$ satisfies a large deviation principle with rate function 
 $
 H(\upsilon|\mu)+H(\upsilon|\nu_0).$
 
 Noting that $\mathcal{W}_{n,\nu_0}$ is the continuous image of $(\mathcal{W}_{n,\nu_0},  \mathcal{L}_{n,\nu_0})$ through projection, we obtain  the result by the contraction principle.
 
  \hfill $\Box$

\begin{remark}\label{3.1}
The Dirichlet posterior given $\eta_1, \ldots, \eta_n$ is a Beta mixture of the Dirichlet process and the conditional distribution of $\{\mathcal{W}_{n,\nu_0}:n\geq 1\}$ given $\xi_i=\eta_i, i=1, \ldots, n$.  Since the coefficient the Dirichlet process part converges to zero as $n$ tends to infinity, the large deviation behaviour of the Dirichlet posterior is the same as the large deviation for the conditional distribution of $\{\mathcal{W}_{n,\nu_0}:n\geq 1\}$. Integrating the posterior with respect to independent product of a probability measure $\nu_0$  leads to an annealed Dirichlet posterior. The large deviation principle in \emph{\cite{ganocon00}} with rate function $H(\nu_0|\mu)$ can be viewed as a quenched large deviation result. In comparison our result in Theorem~\ref{DW} can be viewed as an annealed large deviation principle for the annealed Dirichlet posterior.  The annealed rate function is smaller than the quenched rate function.
 \end{remark}

Recall that the large deviation rate functions for the empirical distributions and the Dirichlet processes are $H(\mu|\nu_0)$ and  $H(\nu_0|\mu)$, respectively.  Since
\[
H(\mu|\nu_0)=H(\mu|\mu)+H(\mu|\nu_0),\ \ H(\nu_0|\mu)=H(\nu_0|\mu)+H(\nu_o|\nu_0),\]
it follows that  $J(\mu,\nu_0)$  is less than both $H(\mu|\nu_0)$ and $H(\nu_0|\mu)$. 

Let 
\beq
{\cal D}_1&=&\{\mu \in M_1(S): H(\mu|\nu_0)<H(\nu_0|\mu)\}\\
{\cal D}_2&=&\{\mu \in M_1(S): H(\mu|\nu_0)>H(\nu_0|\mu)\}\\
{\cal D}_3&=&\{\mu \in M_1(S): H(\mu|\nu_0)=H(\nu_0|\mu)\}.
\eeq
Consider the value of rate function at $\mu$ as the energy required to move from $\nu_0$ to $\mu$. Then it requires less energy to move into region ${\cal D}_1$ for the empirical distributions than the Dirichlet processes. The deviation into ${\cal D}_2$ is easier for  the Dirichlet processes than the empirical distributions. There will be no difference between the two  for region ${\cal D}_3$. The energy required for a deviation of $\{\mathcal{W}_{n,\nu_0}:n\geq 1\}$ from $\nu_0$ is always the lowest among the three.   

For any $\mu, \nu$ in $M_1(S)$, the relative entropy $H(\mu|\nu)$ originated in information theory as the Kullback-Leibler divergence. It describes the information gain when $\mu$ is used in comparison with $\nu$. It is non-negative, convex, and equals to zero only when $\mu=\nu$. But it is not a metric since it is asymmetric, and does not obey  the triangle inequality.  The asymmetry can be rectified by considering the symmetric divergence  $H(\mu|\nu)+H(\nu|\mu)$. But the triangle inequality still does not hold.  However, one does have the following Pinsker's inequality:
 \[
\rho^2_{tv}(\mu,\nu)\leq 2\min\{H(\mu|\nu), H(\nu|\mu)\}
  \]
where $\rho_{tv}(\mu,\nu)$ denote the total variation distance between $\mu$ and $\nu$.  

The function $J(\mu,\nu)$ defined as in \rf{DW-rate} provides a new measurement of divergence between $\mu$ and $\nu$, which we call the $J$-divergence. It is non-negative, symmetric, convex, and equals to zero only when $\mu=\nu$.  The finiteness of  $J(\mu,\nu)$ does not require the absolute continuity between $\mu$ and $\nu$. This can be seen from the following example:
\[
\upsilon(d\,x)=d\,x, \  \  \mu(d\,x)=\frac{1}{2}[d\,x +\delta_{\{0\}}(d\,x)],\  \  \nu(d\,x)=\frac{1}{2}[d\,x +\delta_{\{1\}}(d\,x)].\]

It follows from direct calculation that
\[
 H(\upsilon|\mu)= \log 2=H(\upsilon|\nu), \ \ \ \  J(\mu,\nu)\leq 2\log 2. \]

Clearly  $H(\mu|\nu)=H(\nu|\mu)=\infty$. Thus $J(\mu,\nu)$ can be strictly less than the minimum of $H(\mu|\nu)$ and  $H(\nu|\mu)$. This helps in quantifying the relative information between probabilities that have no absolute continuity relation.  

The $J$-divergence  does not obey the triangle inequality either. But it satisfies an inequality similar to Pinsker's inequality.

\begin{theorem}\label{j-t1}
For any $\mu, \nu$ in $M_1(S)$,
\be\label{j-e1}
\rho^2_{tv}(\mu,\nu) \leq 4 J(\mu,\nu). 
\ee
\end{theorem}
 \proof  One can prove the result using Pinsker's inequality. We choose a direct proof below. 
 
 Fix $\mu, \nu$ in $M_1(S)$. For any $\upsilon$ in $M_1(S)$, and  any measurable subset $A$ of $S$, set
 
 \beq
 &&p=\upsilon(A),\ q=\mu(A),\ \nu(A)= r,\\
 &&\pi(\upsilon)=(p,1-p), \ \pi(\mu)=(q,1-q), \  \pi(\nu)=(r,1-r),
 \eeq
 and 
 \[
 L(p)= H(\pi(\upsilon)|\pi(\mu))+H(\pi(\upsilon)|\pi(\nu)). 
 \]
  
  If $q=r$, then $|\mu(A)-\nu(A)|=0$. If $q=0, 0<r<1$ (the case $0<q<1,r=0$ is similar), then $L(p)$ achieves minimum at $p=0$, and 
  \[
  L(0)= -\log (1-r) \geq r=|\mu(A)-\nu(A)|\geq |\mu(A)-\nu(A)|^2.
  \] 
  
 Next we assume that $0<q,r<1$. 
 
 Solving the equation 
 $L'(p)=0$
 we obtain
 \[
 \hat{p}= \frac{\sqrt{\frac{q r}{(1-q)(1-r)}}}{1+\sqrt{\frac{q r}{(1-q)(1-r)}}}.
 \]
 
Since $L''(p)\geq 0$, it follows that $\hat{p}$ is the minima for $L(\cdot)$, and  
 \[
 L(\hat{p})=-2\log[\sqrt{(1-q)(1-r)}+\sqrt{qr}].
 \]
 
 Without loss of generality we assume $q> r$ and set $\delta=q-r$. It is easily checked by direct calculation that the function
 \[
 \sqrt{(1-q)(1-r)}+\sqrt{qr}=\sqrt{(1-r)^2-\delta (1-r)}+\sqrt{r^2+\delta r} \]
 reaches its maximum $\sqrt{1-\delta^2}$ at $r=\frac{1-\delta}{2}$.  Thus
 \[
  L(\hat{p})\geq -\log[1-\delta^2] \geq \delta^2.\]
  
  Putting all these together we obtain
  \beq
  J(\mu,\nu) &\geq&  \inf\{H(\pi(\upsilon)|\pi(\mu))+H(\pi(\upsilon)|\pi(\nu)): \upsilon \in M_1(S)\}\\
  &\geq& \sup\{|\mu(A)-\nu(A)|^2: A\subset S \}= \frac{\rho^2_{tv}(\mu,\nu)}{4}
  \eeq
  and thus \rf{j-e1}
  
 \hfill $\Box$
 
\section{Applications}

   \setcounter{equation}{0}  
   
   For any $f$ in $C(S)$, consider the random weighted averages $\langle\mathcal{L}_{n,\nu_0},f\rangle$, $\langle \mathcal{V}_{n,\nu_0}, f \rangle$, and $\langle \mathcal{W}_{n,\nu_0}, f \rangle$.  Let $\xi$ have distribution $\nu_0$. It follows from direct calculation that
\beq
&&\mathbb{E}[\langle\mathcal{L}_{n,\nu_0},f\rangle]=\mathbb{E}[\langle \mathcal{V}_{n,\nu_0}, f \rangle] =\mathbb{E}[\langle \mathcal{W}_{n,\nu_0}, f \rangle]=\langle \nu_0, f\rangle,\\
&&\mbox{Var}[\langle\mathcal{L}_{n,\nu_0},f\rangle] = \frac{1}{n}\mbox{Var}[f(\xi)],\\
&&\mbox{Var}[\langle\mathcal{V}_{n,\nu_0},f\rangle] = \frac{1}{n}\mbox{Var}[f(\xi)] +\frac{n-1}{n(n+1)}\langle \nu_0,f^2\rangle\\
&&\mbox{Var}[\langle \mathcal{W}_{n,\nu_0},f\rangle] = \frac{1}{n+1}\mbox{Var}[f(\xi)].
\eeq

Thus all three sequences converge to $\langle \nu_0, f\rangle$ in probability as $n$ tends to infinity. The sequence $\{\langle\mathcal{V}_{n,\nu_0},f\rangle: n\geq 1\}$ has a bigger variance among the three.
 
For  $f(x)=x$, $\langle\mathcal{L}_{n,\nu_0},f\rangle,  \langle \mathcal{V}_{n,\nu_0}, f \rangle$, and $\langle \mathcal{W}_{n,\nu_0}, f \rangle$ are, respectively, the sample mean, the well studied random Dirichlet mean (\cite{Lijoi-Prunster09}), and  the finite Dirichlet weighted mean.  The objective of this section is to compare the large deviation behaviours of the three means when $n$ tends to infinity.  

  Since the map 
\[
\Phi_f: M_1(S) \ra \mathbb{R}, \ \mu \ra \langle \mu, f\rangle
\]
is continuous,  the next result  follows by a direct application of Theorem~\ref{DW}, the large deviations for the Dirichlet process, and the contraction principle.

\begin{theorem}\label{mean-t1}
 Both the family  of laws of $\{\langle \mathcal{V}_{n,\nu_0}, f \rangle:n\geq 1\}$ and the family of laws of $\{\langle \mathcal{W}_{n,\nu_0}, f\rangle:n\geq 1\}$ satisfiy large deviation principles with respective good rate functions
 \[\label{mean-e1}
 I_{1,f}(u; \nu_0)= \inf_{\mu\in M_1(S)}\{H(\nu_0|\mu): \langle\mu, f\rangle=u\}
 \]
 and 
\[\label{mean-e2}
I_{2,f}(u;\nu_0) =\inf_{\mu\in M_1(S)}\{J(\mu,\nu_0): \langle \mu, f\rangle=u\}.
\]

Let $\nu_1$ denote the law of  $M_{1,\nu_0}=\langle\mathcal{V}_{1,\nu_0}, x\rangle$, $I_{i}(\cdot,\nu)=I_{i,f}(\cdot,\nu) $ for $f(x)=x, i=1,2.$ Then we have  
\be\label{mean-e3}
I_{1}(u;\nu_0) =I_{2}(u; \nu_1).
\ee
\end{theorem}

\proof  It suffices to prove \rf{mean-e3}. Let $\{\xi_i: i\geq 1\}$ be i.i.d copies of $M_{1,\nu_0}$ with common distribution $\nu_1$. It follows from \rf{decom}
 that
\beq
M_{n,\nu_0}&= & \langle \mathcal{V}_{n,\nu_0}, x \rangle\\
&=&\frac{1}{\sum_{k=1}^n W_k}\sum_{i=1}^n W_i \xi_i\\
&=& \langle \mathcal{W}_{n,\nu_1}, x\rangle.\eeq

Since both $I_{1,f}(u;\nu_0)$ and $I_{2,f}(u; \nu_1)$ are good rate functions, the result follows from the uniqueness of good rate function. 

\hfill   $\Box$
          
The rate function $I_{1}(u)$ involves the minimization of the relative entropy $H(\nu_0|\mu)$ over the set of probability measures 
$\{\mu \in M_1(S): \langle\mu,x\rangle=u\}$. This is called the reverse information projection in comparison with  the forward information projection of minimizing $H(\mu|\nu_0)$ for $\mu$ over a a subset of $M_1(S)$. The forward information projection over convex set has an explicit solution (\cite{Csiszar75}), and the reverse information projection has a solution over any log-convex domain(\cite{Csiszar03}).  Since the set  $\{\mu \in M_1(S): \langle\mu,x\rangle=u\}$ is not log-convex, 
it is not clear whether one can identify the minimizer, and thus $I_{1}(u)$, explicitly in general.   Our final result will  give an explicit form for $I_{1}(u)$ in the case $\nu_0(d\,x)= d\,x$.

 \begin{theorem}\label{mean-t2}
Let $\nu_0$ be the uniform distribution over $S$. Then 
\beqn\label{mean-e4}
I_1(u,\nu_0)&=& \int_0^1\log (\alpha+\beta x)\nu_0(d\,x) \\
&=& F^{-1}(u)(u-1)+\log (1+F^{-1}(u)u),\nn
\eeqn
where
 \[
\alpha+\beta=\alpha e^{\beta}, \ \ 
 \alpha+\beta u=1.\]
 and 
\[
F(\la) =\left\{\begin{array}{ll}
  \frac{1}{2}&  \la =0,\\
  \frac{e^\la}{e^\la-1} -\frac{1}{\la}& \mbox{ else}.
 \end{array}
 \right.
\]
\end{theorem}      

\proof For any $\mu$ in $M_1(S)$, let $\mu=\mu_1+\mu_2$ be the Lebesgue decomposition  of $\mu$ with respect to $\nu_0$, with 
$\mu_1\ll \nu_0$ and $\mu_2\perp \nu_0$. Set $f(x)=\frac{d\,\mu_1}{d\,\nu_0}$, the Radon-Nikodym derivative of $\mu_1$ with respect to $\nu_0$. Then for any $u$ in $S$ we have  

\beqn\label{mean-e5}
 &&I_{1}(u; \nu_0)=\inf\{\langle\nu_0,-\log f(x)\rangle: \langle \mu_1+\mu_2, x\rangle=u, \mu_1\equiv \nu_0\}\nn\\
 &&=\inf\{\langle \nu_0,-\log f(x)\rangle: f >0, \ \nu_0-a.s., \langle\nu_0, f(x)\rangle \leq 1, \langle \nu_0, x f(x)\rangle \leq u\}\\
 &&=\inf_{a\in [0,1], b\in [0,a\wedge u]}\inf_{ f \in \Gamma_{a,b}}\{\langle\nu_0,-\log f(x)\rangle\} \nn\eeqn
where
\[
\Gamma_{a,b}=\{f \in B(S): f>0, \ \nu_0-a.s.,   \langle \nu_0, f(x)\rangle=a, \langle \nu_0,x f(x)\rangle=b\}.
\]

If $b=0$, then $f(x)=0$ and $\langle \nu_0,-\log f(x)\rangle=\infty$.

 For any $0< a\leq 1$, and $0< b \leq u\wedge a$, let $\la_1(a,b)$ and $\la_2(a,b)$ be the solution to the following equations

\[
\la_1+\la_2=\la_1 e^{a\la_2}, \ \ 
a \la_1+b\la_2=1.
\]

Since $a\geq b$, it follows that
\[
\la_1(a,b)>0,\ \la_1(a,b)+\la_2(a,b)>0.
\]
Thus the nonnegative function 

\[
g_{a,b}(x)=\frac{1}{\la_1(a,b)+\la_2(a,b) x}
\]
is well-defined. It follows from direct calculation that
\[
\langle \nu_0, g_{a,b}(x)\rangle = a, \  \langle \nu_0,x g_{a,b}(x)\rangle=b.
\]
Thus $g_{a,b}$ is in $\Gamma_{a,b}$. For any $f$ in $\Gamma_{a,b}$, we obtain 
\[
\int_0^1 \frac{f(x)}{g_{a,b}(x)}\nu_0(d\,x)= \int_0^1(\la_1(a,b)+\la_2(a,b) x)f(x)\nu_0(d\,x)=1
\]
and 
\[
\int_0^1\log\frac{f(x)}{g_{a,b}(x)}\nu_0(d\,x)\leq\log\bigg(\int_0^1\frac{f(x)}{g_{a,b}(x)}\nu_0(d\,x)\bigg)=0. 
\]
Hence
\[
\langle \nu_0,-\log f(x)\rangle \geq \langle \nu_0, -\log g_{a,b}(x)\rangle,
\]
and  the infimum of $\langle \nu_0,-\log f(x)\rangle$ is achieved at $g_{a,b}$. It is not difficult to see that the map from  $(a,b)$  to $(\la_1,\la_2)$ is one-to-one, and   

\[
\frac{\partial a }{\partial \la_1}<0, \ \ \frac{\partial a}{\partial \la_2}<0, \ \ \frac{\partial b}{\partial \la_1}<0, \ \ \frac{\partial b}{\partial \la_2}<0.\]
It follows  that $\langle \nu_0, -\log g_{a,b}(x)\rangle$ is decreasing in both $a$ and $b$, and the infimum is achieved for $a=1, b=u$. Thus by \rf{mean-e5} we obtain
\[
I_{1}(u; \nu_0)= -\int_{0}^1\log g_{1,u}(x)\nu_0(d\,x)
\]
which implies \rf{mean-e4} and the  theorem by taking $\alpha=\la_1(1,u), \beta=\la_2(1,u)$.

\hfill $\Box$

\begin{remark} \label{4.1}
The function $F(\cdot)$ is the cumulative distribution function of a  random variable.

\end{remark}

\begin{remark}\label{4.2}
Let $\nu_0(d\,x)=d\,x$, and $\{\xi_i:i\geq 1\}$ be i.i.d. with common distribution $\nu_0$. Then the large deviation rate function for the sample mean
\[
\frac{1}{n}\sum_{i=1}^n \xi_i=\langle \mathcal{L}_{n,\nu_0}, x\rangle
\]  
is 
 $$ I_3(u,\nu_0)=\inf\{H(\mu|\nu_0): \langle \mu, x\rangle=u\}.$$ By the minimum discrimination information theorem \emph{(pages 36--39 in \cite{Kullback59})}, the infimum in $I_3(u)$ is achieved at measure $\mu_0$ satisfying
 \[
 \mu_0(d\,x)=c e^{rx}\nu_0(d\,x).
 \]
The constraints 
\[
\langle \nu_0, c e^{rx}\rangle =1,\ \ \langle \nu_0, cx e^{rx}\rangle =u \]
imply that
\[
c=\frac{r}{e^r-1},\  ce^r-1=u r,
\]
and 
\[
u=\frac{e^r}{e^r-1}-\frac{1}{r}=F(r).
\]
By direct calculation,

\beq
I_3(u)&=& H(\mu_0|\nu_0)
= \int_0^1 ce^{rx}(\log c+rx)d\,x\\
&=& \log c + ce^r-1
=r(u-1)+ \log (1+ru)\\
&=&F^{-1}(u)(u-1)+ \log (1+F^{-1}(u)u)\eeq
which is the same as $I_1(u,\nu_0).$ Thus the sample means and the Dirichlet means have the same large deviation behaviour.
\end{remark}

\begin{remark}\label{4.3}
It is known \emph{(\cite{regzi02})} that the law of the Dirichlet mean $M_{n,\nu_0}$ for $n\geq 2$ is absolutely continuous with respect to $\nu_0$ with the Radon-Nikodym derivative
\[
q_n(x)=\frac{n-1}{\pi}\int_{-\infty}^x(x-y)^{n-2}e^{-n\int_0^1\log |y-z|d\,z}\sin(n\pi y)d\,y.   
\]
But it is not clear how one can use this to obtain the explicit form of $I_1(u,\nu_0)$.

\end{remark}

      \end{document}